\documentclass[11pt]{amsart}
\numberwithin{equation}{section}
\newtheorem{theorem}{Theorem}

\numberwithin{theorem}{section} \numberwithin{lemma}{section}
\numberwithin{proposition}{section}

\def\al{\aligned}
\def\eal{\endaligned}

\begin{document}

\tracingpages 1
\title[gradient estimate]{\bf sharp gradient estimate and Yau's Liouville
 theorem for the heat equation on noncompact manifolds}
\author{Philippe Souplet and Qi S. Zhang}
\address{Department of Mathematics, University of Picardie, 02109
St-Quentin,
  France
and Department of Mathematics, University of California, Riverside,
CA 92521, USA }
\date{January 2005}

\begin{abstract}
 We  derive a sharp, localized version of elliptic type gradient estimates
 for positive solutions (bounded or not)
 to the heat equation. These estimates are akin to the Cheng-Yau estimate
for
 the Laplace equation and Hamilton's estimate for bounded solutions to
  the heat equation
 on compact manifolds. As applications, we generalize Yau's celebrated
 Liouville theorem  for positive harmonic functions to
 positive eternal solutions of the heat equation, under certain growth
  condition.  Surprisingly,
 this Liouville theorem for the heat equation
 does not hold even in ${\bf R}^n$ without such a condition.  We also prove
a sharpened long
 time gradient estimate for the log of heat kernel on noncompact
manifolds. This has been an open problem in view of the well known
estimates in the compact, short time case.

\end{abstract}
\maketitle
\section{ {\bf Introduction and main results}}

In recent years, the following three types of gradient estimates
have emerged to play a fundamental role in the study of the Laplace
and the heat equation on manifolds.

\medskip

 {\it {\bf Theorem }(Cheng-Yau [CY]). Let ${\bf M}$ be a complete
manifold with dimension $n \ge 2$, $Ricci ({\bf M}) \ge -k$, $k \ge
0$. Suppose $u$ is any positive harmonic function in  a geodesic
ball $B(x_0, R) \subset {\bf M}$. There holds
\[
 \frac{|\nabla u|}{u} \le  \frac{c_n}{R} + c_n
\sqrt{k}, \leqno(1.1)
\]in $B(x_0, R/2)$, where $c_n$ depends only on the dimension $n$.
}

\medskip

 {\it {\bf Theorem }(Li-Yau [LY]). Let ${\bf M}$ be a complete
manifold with dimension $n \ge 2$, $Ricci ({\bf M}) \ge -k$, $k \ge
0$. Suppose $u$ is any positive solution to the heat equation in
$B(x_0, R) \times [t_0-T, t_0]$. Then
\[
\frac{|\nabla u|^2}{u^2} -\frac{u_t}{u} \le  \frac{c_n}{R^2}+
\frac{c_n}{T}+ c_n k, \leqno(1.2)
\]in $B(x_0, R/2) \times [t_0-T/2, t_0]$. Here $c_n$ depends only on the
dimension $n$. }
\bigskip

{\it {\bf Theorem }(Hamilton [H]). Let ${\bf M}$ be a compact
manifold without boundary and with  $Ricci({\bf M}) \geq -k$, $k\geq
0$ and $u$ be a smooth positive solution of the heat equation with
$u \le M$ for all $(x, t) \in {\bf M} \times (0, \infty)$. Then
\[
 \frac{|\nabla u|^2}{u^2} \le (\frac{1}{t} + 2 k )  \ln \frac{M}{u}.
\leqno(1.3)
\]}
\medskip

Clearly estimate (1.2) reduces to (1.1) when $u$ is independent of
time. On the other hand, for a time dependent solution of the heat
equation, it is well known that the elliptic type gradient estimate
(1.1) can not hold in general. This can be seen from the simple
example where $u(x, t) =  e^{-|x|^2/4t}/(4 \pi t)^{n/2}$ being the
fundamental solution of the heat equation in ${\bf R}^n$. The
parabolic Harnack inequality also exhibits the same phenomenon in
that the temperature at a given point in space time is controlled
from the above by the temperature at a later time. However, in the
case of compact manifolds, Hamilton's estimate (1.3) shows that one
can compare the temperature of two different points at the same time
provided the temperature is bounded.

Just like the Cheng-Yau and Li-Yau estimates, it would be highly
desirable to have a noncompact or localized version of Hamilton's
estimate. However, the example in Remark 1.1 below shows that the
suspected noncompact version of  Hamilton's estimate is {\it false}
even for ${\bf R}^n$ ! This situation contrasts sharply with the
Cheng-Yau and Li-Yau inequality for which the local and noncompact
versions are readily available.

Nevertheless in this paper we discover that, for noncompact
manifolds the elliptic Cheng-Yau estimate actually holds for the
heat equation, after inserting a necessary logarithmic correction
term. This correction term is slightly bigger than that in
Hamilton's theorem (in the power of the log term). However the
estimate holds for noncompact manifolds and it also has a localized
version as the Cheng-Yau estimate. This result seems surprising
since it enables the comparison of temperature distribution
instantaneously, without any lag in time, even for noncompact
manifolds, regardless of the boundary behavior (see Remark 2.1). In
some cases, our estimate (see (1.5) below) even holds for {\it any}
positive solutions, bounded or not. This results seems new even in
${\bf R}^n$ or compact manifolds.

Here is the statement of the theorem.

\begin{theorem}
Let ${\bf M}$ be a Riemannian manifold with dimension $n \ge 2$,
$Ricci ({\bf M}) \ge -k$, $k \ge 0$. Suppose $u$ is any positive
solution to the heat equation in $Q_{R, T} \equiv B(x_0, R) \times
[t_0-T, t_0] \subset {\bf M} \times (-\infty, \infty)$.
 Suppose also $u \le M$ in  $Q_{R, T}$. Then there exists a dimensional
 constant $c$ such that
\[
 \frac{|\nabla u(x, t) |}{u(x, t)} \le c (\frac{1}{R} +
\frac{1}{T^{1/2}}+\sqrt{k})
  \big{(} 1+ \ln \frac{M}{u(x, t)} \big{)}
   \leqno(1.4)
\]in $Q_{R/2, T/2}$.

Moreover, if ${\bf M}$ has nonnegative Ricci curvature and $u$ is
any positive solution of the heat equation on ${\bf M} \times (0,
\infty)$, then there exist dimensional constants $c_1, c_2$ such
that
\[
\frac{|\nabla u(x, t) |}{u(x, t)} \le c_1 \frac{1}{t^{1/2}}
  \big{(} c_2 + \ln \frac{u(x, 2 t)}{u(x, t)} \big{)}
   \leqno(1.5)
\]for all $x \in {\bf M}$ and $t>0$.
\end{theorem}

\medskip

An immediate application of the theorem is the following
time-dependent Liouville theorem, generalizing Yau's celebrated
Liouville theorem for positive harmonic functions, which states that
any positive harmonic function on a noncompact manifold with
nonnegative Ricci curvature is a constant. One tends to expect that
Yau's Liouville theorem would still hold for positive eternal
solutions to the heat equation. However the following simple example
shows that this expectation is false. Let $u=e^{x+t}$ for $x \in{\bf
R}^1.$ Clearly, $u$ is a positive eternal solution for the heat
equation in ${\bf R}^1$ and it is not a constant. Nevertheless, our
next theorem shows that under certain growth conditions, Yau's
Liouville theorem continues to hold for positive eternal solutions
of the heat equation. Moreover, our growth condition in the spatial
direction is sharp by the above example.

\begin{theorem}
 Let ${\bf M}$ be a complete, noncompact manifold with
 nonnegative Ricci curvature. Then the following conclusions hold.

(a).Let $u$ be a positive eternal
 solution to the heat
 equation (i.e. solutions defined in all space-time) such that $u(x,
 t) =  e^{ o (d(x) + \sqrt{|t|})}$ near infinity.
 Then $u$ is a constant.

(b). Let $u$ be an eternal
 solution to the heat
 equation such that $u(x,
 t) = o \big{(} [d(x) + \sqrt{|t|}]
 \big{)}$ near infinity. Then $u$ is a constant.
 \end{theorem}
 \medskip

 Note that the growth condition on the second statement of Theorem 1.2 is
 also sharp in the spatial direction, due to the example $u=x$.
 \medskip

{\it Remark 1.1.} Here we give an example showing that Theorem 1.1
is sharp for noncompact manifolds.  This is suprising since it shows
 that Hamilton's estimate for the compact case (1.3) is actually
false for noncompact manifolds.

For $a>0$ consider $u=e^{a x + a^2 t}$. Clearly $u$ is a positive
solution of the heat equation in $Q=[1, 3] \times [1, 2] \subset
{\bf R} \times (-\infty, \infty).$ Also $\nabla u(2, 2)/u(2, 2) = a$
and $M=\sup_{Q} u = e^{3 a + 2 a^2}$. Hence $\log (M/u(2, 2)) = a$.
Therefore at $(x, t)=(2, 2)$ the left hand side and right hand side
of (1.4) with $R=T=1$ are $a$ and $c(1+a)$ respectively. Obviously
they are equivalent when $a$ is large.

\medskip

Another application of Theorem 1.1 is the following sharpened
estimate on the gradient of the fundamental solution of the heat
equation.

\begin{theorem}
 Let ${\bf M}$ be a $n$ dimensional complete, noncompact manifold with
 nonnegative Ricci curvature. Let $G(x, y, t)$ be the fundamental
 solution of the heat equation. Then there exists  $c=c(n)>0$ such that
\[
\frac{|\nabla_x G(x, y, t) |}{G(x, y, t)} \le c  \frac{1}{t^{1/2}}
  \big{(} 1+
  \frac{d(x, y)^2}{t} \big{)}.
\leqno(1.6)
 \]holds for all $x, y \in {\bf M}$ and $t>0$.
 \end{theorem}
 \medskip

Remark 1.2. For a noncompact manifold and large time, it seems that
all previous results establish an upper bound for $|\nabla G|$
rather than $|\nabla G|/G$. Notice also that estimate (1.6) is
almost sharp since, when ${\bf M}$ is the Euclidean space, one has
\[
\frac{|\nabla_x G(x, y, t) |}{G(x, y, t)} =   \frac{1}{2 \sqrt{t}}
  \frac{d(x, y)}{\sqrt{t}}.
\]We should mention that in the compact and short time case,
estimate similar to (1.5) with a sharp power on the $d(x, y)^2/t$
term have been obtained for all derivatives of $\log G$. See the
papers [S], [H] (Corollary 1.3), [N], [MS], [Hs] and [ST]. The
noncompact case has been wide open.

\medskip

Our work is motivated by [CY], [LY], [Y1-2], [H] and [BL] and [GGK].
\medskip

\section{\bf Proof of Theorems 1.1, 1.2 and 1.3}


The order of proofs is Theorem 1.1, 1.3 and 1.2.
\medskip

{\bf Proof of Theorem 1.1.}

Suppose $u$ is a solution to the heat equation in the statement of
the theorem in the parabolic cube $Q_{R, T}=B(x_0, R) \times [t_0-T,
t_0]$. It is clear that the gradient estimate in Theorem 1.1 is
invariant under the scaling $u \to u/M$. Therefore, we can and do
assume that $0< u \le 1$.

Write
\[
f = \ln u, \qquad w \equiv | \nabla \ln (1-f)|^2 = \frac{|\nabla
f|^2}{(1-f)^2}. \leqno(2.1)
\]

Since $u$ is a solution to the heat equation, simple calculation
shows that
\[
\Delta f + | \nabla f|^2 - f_t =0, \leqno(2.2)
\]

We will derive an equation for $w$. First notice that
\[
\al w_t &= \frac{2 \nabla f (\nabla f)_t}{(1-f)^2} + \frac{2
|\nabla f|^2 f_t}{(1-f)^3}\\
&= \frac{2 \nabla f \nabla (\Delta f + |\nabla f|^2)}{(1-f)^2} +
\frac{2 |\nabla f|^2 (\Delta f + |\nabla f|^2)}{(1-f)^3} \eal
\]In local orthonormal system, this can be written as
\[
w_t=\frac{2  f_j f_{iij} + 4 f_i f_j f_{ij}}{(1-f)^2} + 2 \frac{
f^2_i f_{jj} + |\nabla f|^4}{(1-f)^3}. \leqno(2.3)
\]Here and below, we have adopted the convention $u^2_i = |\nabla
u|^2$ and $u_{ii} = \Delta u$.

Next
\[
\nabla w = \big{(} \frac{f^2_i}{(1-f)^2} \big{)}_j = \frac{2 f_i
f_{ij}}{(1-f)^2} + 2 \frac{f^2_i f_j}{(1-f)^3}. \leqno(2.4)
\]It follows that
\[
\al \Delta w &= \big{(} \frac{f^2_i}{(1-f)^2} \big{)}_{jj}\\
&= \frac{ 2 f^2_{ij}}{(1-f)^2} + \frac{2 f_i f_{ijj}}{(1-f)^2} +
\frac{4 f_i f_{ij} f_j}{(1-f)^3}\\
&\qquad + \frac{4 f_i f_{ij} f_j}{(1-f)^3} + 2 \frac{f^2_i
f_{jj}}{(1-f)^3} + 6 \frac{f^2_i f^2_j}{(1-f)^4}. \eal \leqno(2.5)
\]By (2.5) and (2.3),
\[
\al
&\Delta w - w_t \\
&=\frac{ 2 f^2_{ij}}{(1-f)^2} + 2 \frac{ f_i f_{ijj}-f_j
f_{iij}}{(1-f)^2}\\
&\qquad + 6 \frac{|\nabla f|^4}{(1-f)^4} + 8 \frac{ f_i f_{ij}
f_j}{(1-f)^3} + 2 \frac{f^2_i f_{jj}}{(1-f)^3} \\
&\qquad - 4 \frac{ f_i f_{ij} f_j}{(1-f)^2} - 2 \frac{f^2_i
f_{jj}}{(1-f)^3} - 2 \frac{|\nabla f|^4}{(1-f)^3}. \eal
\]The 5th and 7th terms on the righthand side of this identity
cancel each other. Also, by Bochner's identity
\[
f_i f_{ijj}-f_j f_{iij} = f_j( f_{jii}-f_{iij}) = R_{ij} f_i f_j \ge
-k |\nabla f|^2,
\]where $R_{ij}$ is the Ricci curvature. Therefore
\[
\al
&\Delta w - w_t \\
& \ge \frac{ 2 f^2_{ij}}{(1-f)^2} + 6 \frac{|\nabla f|^4}{(1-f)^4} +
8 \frac{ f_i f_{ij} f_j}{(1-f)^3} - 4 \frac{ f_i f_{ij}
f_j}{(1-f)^2} - 2 \frac{|\nabla f|^4}{(1-f)^3} -\frac{2 k |\nabla
f|^2}{(1-f)^2}. \eal \leqno(2.6)
\]Notice from (2.4) that
\[
\nabla f \nabla w = \frac{2 f_i f_{ij}f_j }{(1-f)^2} + 2 \frac{f^2_i
f^2_j}{(1-f)^3}.
\]Hence
\[
0 = 4 \frac{ f_i f_{ij} f_j}{(1-f)^2}  - 2 \nabla f \nabla w + 4
\frac{|\nabla f|^4}{(1-f)^3}, \leqno(2.7)
\]
\[
0= -4 \frac{ f_i f_{ij} f_j}{(1-f)^3} + [ 2 \nabla f \nabla w - 4
\frac{|\nabla f|^4}{(1-f)^3} ] \frac{1}{1-f}. \leqno(2.8)
\]Adding (2.6) with (2.7) and (2.8), we deduce
\[
\al
&\Delta w - w_t \\
& \ge \frac{ 2 f^2_{ij}}{(1-f)^2} + 2 \frac{|\nabla f|^4}{(1-f)^4}
+ 4 \frac{ f_i f_{ij} f_j}{(1-f)^3}  \\
&\qquad + \frac{2}{1-f} \nabla f \nabla w - 2 \nabla f \nabla w + 2
\frac{|\nabla f|^4}{(1-f)^3} -\frac{2 k |\nabla f|^2}{(1-f)^2}. \eal
\]Since
\[
\frac{ 2 f^2_{ij}}{(1-f)^2} + 2 \frac{|\nabla f|^4}{(1-f)^4} + 4
\frac{ f_i f_{ij} f_j}{(1-f)^3} \ge 0,
\]we have
\[
\Delta w - w_t
 \ge  \frac{2 f}{1-f} \nabla f \nabla w  +
2 \frac{|\nabla f|^4}{(1-f)^3} -\frac{2 k |\nabla f|^2}{(1-f)^2}.
\]Since $f \le 0$, it follows that
\[
\Delta w - w_t
 \ge  \frac{2 f}{1-f} \nabla f \nabla w  +
2 (1-f) \frac{|\nabla f|^4}{(1-f)^4}-\frac{2 k |\nabla
f|^2}{(1-f)^2},
\]i.e
\[
\Delta w - w_t
 \ge  \frac{2 f}{1-f} \nabla f \nabla w  +
2 (1-f) w^2 - 2k w. \leqno(2.9)
\]From here, we will use the well known cut-off function by Li-Yau [LY],
to derive the desired bounds. We caution the reader that the
calculation is not the same as in that of [LY] due to the difference
of the first order term.

 Let $\psi=\psi(x, t)$ be a smooth cut-off function
supported in $Q_{R, T}$, satisfying the following properties

(1). $\psi = \psi(d(x, x_0), t) \equiv \psi(r, t)$; $\psi(x, t) = 1$
in $Q_{R/2, T/4}$, \ $0 \le \psi \le 1$.

(2). $\psi$ is decreasing as a radial function in the spatial
variables.

(3). $\frac{|\partial_r \psi|}{\psi^{a}} \le \frac{C_a}{R}$,
$\frac{|
\partial^2_r \psi|}{\psi^{a}} \le \frac{C_a}{R^2}$when $0<a<1$.

(4). $\frac{|\partial_t \psi|}{\psi^{1/2}} \le \frac{C}{T}$.

 Then, from
(2.9) and a straight forward calculation, one has
\[
\aligned
 \Delta (\psi w)& + b \cdot \nabla (\psi w) - 2 \frac{\nabla
 \psi}{\psi}
 \cdot  \nabla (\psi w) - (\psi w)_t\\
&\ge 2 \psi (1-f) w^2 + (b  \cdot \nabla \psi) w - 2 \frac{|\nabla
\psi|^2}{\psi} w + (\Delta \psi) w -\psi_t w - 2k w \psi,
\endaligned
\leqno(2.10)
\]where we have written
\[
b = - \frac{2 f}{1-f} \nabla f.
\]

Suppose the  maximum of $\psi w$ is reached at $(x_1, t_1)$. By
[LY], we can assume, without loss of generality that $x_1$ is not in
the cut-locus of ${\bf M}$. Then at this point, one has,
 $\Delta (\psi w) \le 0$, $(\psi w)_t \ge 0$ and $\nabla (\psi
 w)=0$. Therefore
 \[
2 \psi (1-f) w^2(x_1, t_1) \le - [ \ (b  \cdot \nabla \psi) w - 2
\frac{|\nabla \psi|^2}{\psi} w + (\Delta \psi) w -\psi_t w \ + 2k w
\psi](x_1, t_1). \leqno(2.11)
\]We need to find an upper bound for each term of the righthand
side of (2.11).

\[
\al |(b  \cdot \nabla \psi) w| &\le \frac{ 2|f|}{1-f} |\nabla f| w
|\nabla \psi| \le 2 w^{3/2} |f| \  |\nabla \psi| \\
&= 2 [\psi (1-f) w^2]^{3/4} \ \frac{f |\nabla \psi|}{[\psi
(1-f)]^{3/4}}\\
&\le \psi (1-f) w^2 + c \frac{(f |\nabla \psi|)^4}{[\psi (1-f)]^3}.
\eal
\]This implies
\[
|(b  \cdot \nabla \psi) w| \le (1-f) \psi w^2 + c \frac{f^4}{R^4
(1-f)^3}. \leqno(2.12)
\]

For the second term on the righthand side of (2.11), we proceed as
follows
\[
\aligned
 \frac{|\nabla \psi|^2}{\psi} w &= \psi^{1/2} w
\frac{|\nabla
\psi|^2}{\psi^{3/2}}\\
&\le \frac{1}{8} \psi w^2 + c \big{(} \frac{|\nabla
\psi|^2}{\psi^{3/2}} \big{)}^2 \le \frac{1}{8} \psi w^2  + c
\frac{1}{R^4}.
\endaligned
\leqno(2.13)
\]

Furthermore, by the properties of $\psi$ and the assumption of on
the Ricci curvature, one has
\[
\aligned - (\Delta \psi) w &= -(\partial^2_r \psi + (n-1) \frac{
\partial_r \psi}{r} + \partial_r \psi \partial_r \ln \sqrt{g} ) w\\
&\le (|\partial^2_r \psi| + 2 (n-1) \frac{|
\partial_r \psi|}{R} + \sqrt k |\partial_r \psi|  ) w\\
 &\le \psi^{1/2} w \frac{|\partial^2_r
\psi|}{\psi^{1/2}} + \psi^{1/2} w 2 (n-1) \frac{|
\partial_r \psi|}{R \psi^{1/2}} +
 \psi^{1/2} w \frac{\sqrt k |\partial_r \psi|}{\psi^{1/2}}\\
&\le \frac{1}{8} \psi w^2  + c \big{(} [\frac{|\partial^2_r
\psi|}{\psi^{1/2}}]^2 + [ \frac{|
\partial_r \psi|}{R \psi^{1/2}}]^2 +
 [ \frac{\sqrt k |\partial_r \psi|}{\psi^{1/2}}]^2 \big{)}.
 \endaligned
\]Therefore
\[
- (\Delta \psi) w \le \frac{1}{8} \psi w^2  + c \frac{1}{R^4} + c k
\frac{1}{R^2}. \leqno(2.14)
\]

Now we estimate $|\psi_t| \  w$.
\[
\aligned
 |\psi_t| \ w &= \psi^{1/2} w \frac{|\psi_t|}{\psi^{1/2}}\\
&\le \frac{1}{8} \big{(} \psi^{1/2} w \big{)}^{2} + c \big{(}
\frac{|\psi_t|}{\psi^{1/2}} \big{)}^{2}.
\endaligned
\]This shows
\[
|\psi_t| w  \le \frac{1}{8} \psi w^2 + c \frac{1}{T^2}. \leqno(2.15)
\]Finally, for the last term, we have
\[
2 k w \psi \le \frac{1}{8} \psi w^2 + c k^2. \leqno(2.16)
\]

Substituting (2.12)-(2.16) to the righthand side of (2.11), we
deduce,
\[
 2 (1-f) \psi w^2 \le (1-f)  \psi w^2 + c \frac{f^4}{R^4
(1-f)^3} + \frac{1}{2} \psi w^2 + \frac{c}{R^4} + \frac{c}{T^2} +
\frac{c k}{R^2} + c k^2.
\]Recall that $f \le 0$, therefore the above implies
\[
 \psi w^2(x_1, t_1) \le  c \frac{f^4}{R^4 (1-f)^4} +
  \frac{1}{2} \psi w^2(x_1, t_1) +
\frac{c}{R^4} +  \frac{c}{T^2} + c k^2.
\]Since $\frac{f^4}{(1-f)^4} \le 1$, the above shows, for all $(x,
t)$ in $Q_{R, T}$,
\[
\al \psi^2(x, t) w^2(x, t) &\le \psi^2(x_1, t_1) w^2(x_1, t_1)
\\
&\le \psi(x_1, t_1) w^2(x_1, t_1) \\
&\le c \frac{c}{R^4} + \frac{c}{T^2} +  c k^2. \eal
\]Notice that $\psi(x, t) = 1$ in $Q_{R/2, T/4}$ and $w =|\nabla
f|^2/(1-f)^2$. We finally have
\[
\frac{|\nabla f(x, t)|}{1-f(x, t)} \le
 \frac{c}{R} + \frac{c}{\sqrt{T}} + c \sqrt{k}.
\]We have completed the proof of (1.4) since $f=\ln (u/M)$ with $M$
scaled to $1$.

To prove (1.5), we apply (1.4) on the cube $Q_{\sqrt{t}, t/2} = B(x,
\sqrt{t}) \times [ t/2, t]$. By Li-Yau's inequality [LY], we know
that
\[
M = \sup_{Q_{\sqrt{t}, t/2}} u \le c u(x, 2t).
\]Now (1.5) follows from (1.4).
\qed

\bigskip

{\bf Proof of Theorem 1.3.}

Let $G(x, y, t)$ be the fundamental solution of the heat equation on
${\bf M}$. By the Li-Yau estimate [LY], given any $\delta>0$, there
exist $c_1, c_2>0$ such that,
\[
\frac{c_2}{|B(x, \sqrt{t})|} e^{- d(x, y)^2/(4-\delta) t} \le G(x,
y, t) \le \frac{c_1}{|B(x, \sqrt{t})|} e^{- d(x, y)^2/(4+\delta) t},
\leqno(2.16)
\]for all $x, y \in {\bf M}$ and $t>0$.
Fixing $x, y$ and $t$, we apply Theorem 1.1 on the function $u(z,
\tau) = G(z, y, \tau)$ on the cube $Q=B(x, \sqrt{t}) \times [t/2,
t]$. By the upper bound in (2.16) and the volume doubling property
of the manifold, we know that
\[
 u(z, \tau) \le \frac{c_1}{|B(z,
\sqrt{\tau})|} \le \frac{c_3}{|B(x, \sqrt{t})|}
\]for $(z, \tau) \in Q$. By Theorem 1.1 and the lower bound in (2.16), we
deduce
\[
 \frac{|\nabla u(x, t) |}{u(x, t)} \le c  \frac{1}{t^{1/2}}
  \big{(} 1+
  \ln \frac{ \frac{c_3}{|B(x, \sqrt{t})|}}{\frac{c_2}{|B(x, \sqrt{t})|}
  e^{- d(x, y)^2/(4-\delta) t}} \big{)}.
 \]That is
\[
\frac{|\nabla_x G(x, y, t) |}{G(x, y, t)} \le c_4  \frac{1}{t^{1/2}}
  \big{(} 1+
  \frac{d(x, y)^2}{t} \big{)}.
 \]\qed
\medskip

{\bf Proof of Theorem 1.2.}

(a).  The proof is immediate. By our assumption the function $u+1$
satisfies $\ln (u+1) = o(d(x)+ \sqrt{t})$ near infinity. Fixing
$(x_0, t_0)$ in space time and using Theorem 1.1 for $u+1$ on the
cube $B(x_0, R) \times [t_0-R^2, t_0]$, we have
\[
\frac{|\nabla u(x_0, t_0)|}{u(x_0, t_0)+1} \le \frac{C}{R} [1 + o(R)
].
\]

Letting $R\to \infty$, it follows that $|\nabla u(x_0, t_0)|=0$.
Since $(x_0, t_0)$ is arbitrary, one sees that $u=c$.

\medskip

(b). Fix $(x_0, t_0)$ in space time and let $A_R = \sup_{Q_{R,
\sqrt{R}}} |u|.$ Consider the function $U = u + 2 A_{2R}$. Clearly
$A_{2R} \le U(x, t) \le 3 A_{2R}$ when $(x, t) \in Q_{2R,
\sqrt{2R}}$.  Applying Theorem 1.1 (a) for $U$, we deduce
\[
\frac{|\nabla u(x_0, t_0)|}{u(x_0, t_0)+2 A_{2R}} \le  \frac{C}{R}.
\]Since $A_{2R} = o(R)$ by assumption, the result follows after
taking $R$ to $\infty$. \qed
\medskip

{\it Remark 2.1.} In case ${\bf M}$ has nonnegative Ricci curvature,
(1.4) reduces to
\[
 \frac{|\nabla u(x, t) |}{u(x, t)} \le c (\frac{1}{R} +
\frac{1}{T^{1/2}})
  \big{(} 1+ \ln \frac{M}{u(x, t)} \big{)}
 \]in $Q_{R/2, T/2}$. If $d(x, y) \le \sqrt{T}$, we can integrate
 this along a geodesic to deduce
 \[
 \frac{M}{u(y, t)} \ge c \big{(} \frac{M}{u(x, t)} \big{)}^{\theta},
 \]for some $\theta \in (0, 1)$ and $c>0$ depending only on $n$.
Note in the compact case, Hamilton's theorem implies the above with
$\theta =1$, which gives an elliptic  Harnack inequality for the
heat equation for compact manifolds. However, as pointed in the
introduction, the elliptic Harnack inequality is false for the heat
equation on noncompact manifolds in general. This explains the
difference between our noncompact gradient estimate (1.4) and
Hamilton's compact estimate (1.3).

\bigskip

{\bf Acknowledgement.} We thank Professors Huaidong Cao, Bennet
Chow, Alexander Grigoryan, Peter Li, Zhiqin Lu and Lei Ni for very
helpful conversations.

\bigskip

\noindent e-mails: souplet@math.uvsq.fr and qizhang@math.ucr.edu

\enddocument